\documentclass[10pt]{article}
\usepackage{latexsym,amsmath,amssymb,mathrsfs,epsfig,graphicx,rotating}
\usepackage[latin1]{inputenc}
\newtheorem{thm}{Theorem}[section]

\newtheorem{defi}{Definition}[section]

\newcommand{\Real}{\bf R}

\begin{document}

\title{\large \bf On the approximation of  derivatives using  divided difference operators \\ [0.3em]
                    preserving the local convergence order of iterative methods}
\author{{Miquel Grau--S\'{a}nchez$^{\rm \,a}$, Miquel Noguera$^{\rm \,a}$, Sergio Amat$^{\rm \,b}$} \\             \hfill    \\ [-0.4em]
 {\footnotesize $\,^{\rm a}${\it Technical University of Catalonia, Department of Applied Mathematics II}}\\
      {\footnotesize {\it Jordi Girona 1-3, Omega, 08034 Barcelona, Spain}} \\
     \hfill    \\ [-0.8em]
     {\footnotesize $\,^{\rm b}${\it Universidad Politécnica de Cartagena, Departamento de Matemática Aplicada y Estadística}} \\
      {\footnotesize {\it  Paseo de Alfonso XIII, 52, 30203 Cartagena (Murcia), Spain}} \\ [0.1em]
      {\footnotesize  E-mail address: $\{$miquel.grau, miquel.noguera$\}$@upc.edu, sergio.amat@upct.es}
      }
\date{}
\maketitle

%
%
%

\vspace{-4mm}
\begin{abstract}
\noindent
A development of an inverse first-order divided difference operator for functions of several variables is presented.
Two generalized derivative-free algorithms builded up from Ostrowski's method for solving systems of nonlinear
equations are written and analyzed. A direct computation of the local order of convergence for these variants of Ostrowski's method is given. In order to preserve the local order of convergence, any divided difference operator
is not valid.
Two counterexamples of computation of a classical divided difference operator without preserving the order
are presented. A new divided difference operator solving this problem is proposed.
Furthermore, a computation that approximates the order of convergence is generated for the examples and it confirms
in a numerical way that the order of the methods is well deduced.
\end{abstract}

%
%

\vspace{1mm}
{\small

\vspace{0.5mm}
}


\section{Introduction}

In this paper, we consider iterative methods to find a simple root of a system of nonlinear equations

\vspace{-2mm}
\begin{equation}\label{equ1}
    F(x) \,= \, 0 ,
\end{equation}
where $\,F:D \subset \Real^m \longrightarrow \Real^m\,$ and $D$ is an open convex domain in $\Real^{m}$.
We assume that the solution of (\ref{equ1}) is $\,\alpha \in D\,$ at which $\, F'(\alpha) \,$ is nonsingular.

\vspace{3mm}
The goal of this paper is to give a generalization for several variables of a given Newton's type
iterative method using an approximation of the derivative by means of a central difference
and preserving the local order of convergence. We center our analysis on Ostrowski's method (\ref{eq1})
and the improvement (\ref{eq2}), but similar analysis should be performed on other methods.

\vspace{3mm}
We would like to mention that the matrix of the divided difference operator, used generally in practice \cite{PoPt,Argy,HeRu,AB2,AB1,Ostro_A},
does not preserve the local order of convergence in general. We present two counterexamples of this fact and
we propose a new divided difference operator preserving the local order.

\vspace{3mm}
After setting our iterative method we recall Ostrowski's iterative function \cite{Ostr} for a single
nonlinear equation, $f(x) = 0$, with fourth-order that is given by

\vspace{-7mm}
\begin{eqnarray} \nonumber
y & = & x -\, \frac{f(x)}{f^{\prime}(x)} , \\ \label{eq1}
X & = & y - \mu \, f(y) ,
\end{eqnarray}
where $f: \Real \longrightarrow \Real\,$,
$ \mu \,=\,1 / \left(2\, [y, x; f] - f^{\prime}(x) \right))\,$ and $\, [y, x; f] = (f(y)-f(x))/(y - x)$ is the
divided difference of $f$ at the points $x$ and $y$.

\vspace{2mm}
A variant of Ostrowski's method with sixth-order convergence is proposed in \cite{GrDi1}. This improvement has
the following scheme

\vspace{-7mm}
\begin{eqnarray} \nonumber
y & = & x -\, \frac{f(x)}{f^{\prime}(x)} , \\ \label{eq2}
z & = & y - \mu \, f(y) , \\ [0.3em] \nonumber
X & = & z - \mu \, f(z) .
\end{eqnarray}
This method improves the local order of convergence of Ostrowski's method with an additional evaluation of
the function. Other works related with the improvement of the local order of convergence of (\ref{eq1}) can
be found in \cite{KLW, BRW,BWR,KWL,KW}.

\vspace{2mm}
In \cite{CHMT} a well-known technique (see \cite{Trau}) that consists in replacing the first derivative
of (\ref{eq1}) and (\ref{eq2}) by a central difference is used in order to obtain a new method that preserves
the fourth and the sixth convergence order respectively. Notice that these iterative methods are derivative
free, since $f^{\prime} (x)$ is approximated by $ \left[\, x+f(x), x-f(x); f \,\right]$.

\vspace{2mm}
The paper is organized as follows:
We start in Section 2 with a generalization for several variables of the Ostrowski's method (\ref{eq1})
and the improvement (\ref{eq2}) using an approximation of the derivative by means of a central difference
and preserving the local order of convergence.
Section 3 is devoted to the computational efficiency index performing a comparison for the three iterative methods.
By means of two counterexamples, in Section 4, we emphasize the necessity of a new way to compute the matrix of divided difference operator since in many examples the theoretical one and the classical way used until now to calculate it do not coincide. We finalize the paper with some concluding remarks in Section 5.


\section{Local order of convergence}

\vspace{1mm}
We assume that $F$ has, at least, fourth-order derivatives with continuity on $D$. We consider the first divided difference operator of $F$ on $\Real^m$ as a mapping   $\;  \left[ \, \cdot ,
\cdot \, ; F \right]  : D \times D \subset \Real^m \times \Real^m \longrightarrow \mathscr{L}(\Real^m)$ which
is defined by (see \cite{OrRh,Argy,GGN} and the references therein)

\vspace{-3mm}
\begin{equation} \label{equivDD}
     \left[\, x+h , x ; F \right]  \,  =  \,  \int_0^1 \: F ' (x + t h) \, dt , \qquad \forall (x, h) \in
     \Real^m \times \Real^m .
\end{equation}

\vspace{-1mm}
Developing  $\,F ' ( x + t\,h ) $ in Taylor's series at the point $x$ and integrating we obtain

\vspace{-5mm}
\begin{eqnarray} \label{Int}
\left[\, x+h , x ; F \right] &=& F^{\prime}(x) \,+\,  \frac{1}{2} \; F^{\prime\prime}(x)\: h
+\,  \frac{1}{6} \; F^{\prime\prime\prime}(x)\: h^2 +\, O(h^3) .
\end{eqnarray}


\vspace{-1mm}
Taking into account that $e = x - \alpha$ we develop $F(x)$ and its derivatives in a neighborhood of $\alpha$.
Namely, assuming that $ \left[ F ' (\alpha) \right]^{-1}$ exists, we have

\vspace{-3.5mm}
\begin{equation}\label{FTay}
    F (x) = \Gamma  \left(  e\, + \,  A_2 \,e^2 + \,  A_3 \,e^3  +\, O(e^4) \right) ,
\end{equation}

where $\, \Gamma =  F ' (\alpha)$ and $\,{\displaystyle A_p = \,\frac{1}{p\,!}\; \Gamma^{-1} F^{(p)}(\alpha) \in
\mathscr{L}_p (\Real^m, \Real^m)}$, $p = 2, \,3$. From (\ref{FTay}) the derivatives of $F(x)$
can be written as

\vspace{-4mm}
\begin{equation} \label{Fs10}
F'(x) =  \Gamma  \left( I  +  2\,A_2 \,e +  3\,A_3 \,e^2 +\, O(e^3)  \right), \; \; \;
F ''(x) = \Gamma  \left( 2\,A_2 +  6\,A_3 \,e +\, O(e^2) \right), \; \; \; \,
F '''(x) = \Gamma  \left( 6\,A_3 +\, O(e) \right).
\end{equation}
Setting $y = x + h $ and $E= y -\alpha$, we have $ h = E - e$.
Replacing the previous expressions (\ref{Fs10}) into  (\ref{Int}) we get

\vspace{-6mm}
\begin{eqnarray} \label{eE11}
\left[ \, y, x ; F \right] \,=\, F ' (\alpha) \left( I +\, A_2 \,( E + e ) +\, A_3 \, (E^2 + E e + e^2 )
+\, O_3(E, e) \right) ,
\end{eqnarray}

\vspace{-1mm}
where $\,  O_3 (E, e)$ represents the terms of order higher than three. We say that a function depending on
$E$ and $ e$ is an $O_3 (E, e)$ if it is an $ O(E^{q_0}\,e^{q_1})$ with $q_0 + q_1 = 3, \; q_i \ge 0, \: i=0,1$.

\vspace{2mm}
In all our analysis we must consider the operator $ \left[ \, x+F(x), x-F(x) ; F \right] $.
From $x+F(x)-\alpha = e + F(x)$ and $x-F(x)-\alpha = e - F(x)$, if we replace in (\ref{eE11}) $E$
by $e + F(x)$, $e$ by $ e - F(x)$, and $F(x)$ by (\ref{FTay}), then we will get

\vspace{-4mm}
\begin{eqnarray} \label{CDD}
\left[ \, x+F(x), x-F(x) ; F \right] \,=\, \Gamma \left( I + 2\, A_2 \, e  +\, A_3 \, ( 3\, e^2 + {\widetilde{e}}^{\,2} ) +\, O(e^3) \right) ,
\end{eqnarray}
where $\widetilde{e} = F ' (\alpha)\, e$. Expanding in formal power developments the inverse of the preceding
operator (\ref{CDD}) in terms of $ e$ and $ \widetilde{e}$  we get

\vspace{-6mm}
\begin{eqnarray} \label{CDDI}
\left[ \, x+F(x), x-F(x) ; F \right]^{-1} \,=\,  \left( I - 2\, A_2 \,e -\,A_3 \,( 3\,e^2 + {\widetilde{e}}^{\,2} ) +\,4\,A^2_2\,e^2 +\, O(e^3) \right) \, \Gamma^{-1} .
\end{eqnarray}

\vspace{2mm}
The first step of the variant of  Ostrowski's method is  Newton's method using a central divided difference.
That is, we have the following iteration function

\vspace{-5mm}
\begin{eqnarray} \label{VN}
y \,=\, \Phi_0 (x) \,=\, x -\, \left[ \, x+F(x), x-F(x) ; F \right]^{-1} \,F(x) .
\end{eqnarray}
Subtracting the root $\alpha$ from both sides of (\ref{VN}) and considering (\ref{FTay}) and (\ref{CDDI}), the corresponding vectorial error equation is given by

\vspace{-7mm}
\begin{eqnarray*}
E\,=\,y - \alpha &=& e -\,\left( I - 2\, A_2 \, e + O(e^2) \right) \left( e\, +\,A_2\,e^2 +\,O(e^3)\right) \\
 &=& A_2 \, e^2 +\, O(e^3) .
\end{eqnarray*}
The second step of the variant of  Ostrowski's method is defined by the following iterative function

\vspace{-5mm}
\begin{eqnarray} \label{VO}
z \,=\, \Phi_1 (x,y)   \,=\, y -\, \nu \,F(y) ,
\end{eqnarray}
where $ \nu = \left( 2\,\left[ \, y, x ; F \right] -\,\left[ \, x+F(x), x-F(x) ; F \right] \right)^{-1}$.
Note that from the developments given in  (\ref{eE11}) and (\ref{CDD}) we have

\vspace{-7mm}
\begin{eqnarray*}
\nu^{-1} &=& \Gamma \left( I + 2\,A_2\,E +\,A_3 \,( 2\,E^2 +\,2\,e\,E -e^2
-\widetilde{e}^{\,2} ) + O(e^3) \right) \\
  &=& \Gamma  \left( I + 2\,A^2_2\,e^2 -\,A_3 \,( e^2 + \widetilde{e}^{\,2} ) + O(e^3) \right) ,
\end{eqnarray*}
and applying again the computation of the inverse operator we have in this case

\vspace{-5mm}
\begin{eqnarray*}
\nu &=& \left( I - 2\,A^2_2\,e^2 +\,A_3 \,( e^2 + \widetilde{e}^{\,2} ) + O(e^3) \right) \, \Gamma^{-1} .
\end{eqnarray*}

Subtracting the root $\alpha$ from both sides of (\ref{VO}) the vectorial error equation for the iterative method defined in (\ref{VN}) and (\ref{VO}) is

\vspace{-7mm}
\begin{eqnarray} \nonumber
\varepsilon\,=\,z - \alpha &=&E -\,\left( I - 2\,A^2_2\,e^2 +\,A_3 \,( e^2 + \widetilde{e}^{\,2} ) + O(e^3)\right) \left( E\, +\,A_2\,E^2 +\,O(E^3)\right) \\ \label{EVO}
 &=& \left( A^2_2 -\,A_3\,A_2 \right) \, e^4 -\, A_3\,\widetilde{e}^{\,2} \,A_2 \,e^2  +\, O(e^5) .
\end{eqnarray}

Summarizing the preceding results we are in a condition to state the following theorem.

\begin{thm}\label{th21}
The iterative method given by the iteration function $\Phi_1$ defined in (\ref{VN}) and (\ref{VO}) is at least of fourth local convergence order and its vectorial error difference equation is given in (\ref{EVO}).
\end{thm}

\vspace{1mm}
Moreover, if we introduce a new step with an additions evaluation of the function $F$ in the following way

\vspace{-6mm}
\begin{eqnarray} \label{VO2}
X \,=\, \Phi_2 (x,y,z)   \,=\,  z -\, \nu \,F(z) ,
\end{eqnarray}

\vspace{-3mm}
then the error becomes

\vspace{-5mm}
\begin{eqnarray} \nonumber
 X- \alpha &=& \varepsilon -\,\left( I - 2\,A^2_2\,e^2 +\,A_3 \,( e^2 + \widetilde{e}^{\,2} ) + O(e^3) \right)
 \left( \varepsilon +\,O(\varepsilon^2)\right) \\ \label{EVO2}
 &=& \left[ \left( 2\,A^2_2 -\,A_3 \right) e^2 -\,A_2\, \widetilde{e}^{\,2} \right] \: \left[ \left( A^2_2
 -\,A_3\,A_2 \right)  e^4 -\, A_3\,\widetilde{e}^{\,2} \,A_2 \,e^2 \right]  +\, O(e^7) ,
\end{eqnarray}

and the local error of convergence of the iterative method defined by the steps (\ref{VN}), (\ref{VO}) and
(\ref{VO2}) is at least of sixth order. We present the preceding results in the following theorem.

\begin{thm}\label{th22}
The iterative method given by the iteration function $\Phi_2$ defined by (\ref{VN}), (\ref{VO}) and (\ref{VO2})
is at least of sixth local order of convergence and the vectorial error difference equation is (\ref{EVO2}).
\end{thm}

%
%
%
%

\section{Computational efficiency and comparison between the methods $\Phi_i$}

The computational efficiency index ($CEI$) and the computational cost per iteration ($ {\cal{C}}$) are defined as

\vspace{-2mm}
\begin{equation}\label{CEIC}
CEI(\mu, m, \ell) \,=\, \rho^{ \frac{1}{{\cal{C}}(\mu, m, \ell)}} , \quad \mbox{and} \quad
{\cal{C}}(\mu, m, \ell) = a(m)\, \mu + p(m, \ell) ,
\end{equation}

\vspace{1mm}
where $a(m)$  represents the number of evaluations of the scalar functions $(F_1, \ldots, F_m)$
used in the evaluation of $F$, $[x-F(x), x+F(x); F]$ and $[y, x; F]$. The number of products needed per iteration is represented by $p(m, \ell)$.
In order to express in terms of products the value of ${\cal{C}}(\mu, m, \ell)$ a ratio $\mu > 0 \,$ between products and scalar functions evaluations and a ratio $\,\ell\ge 1\,$ between  products and quotients are required (see \cite{GGN,GGN2}). A first and classical definition of difference divided operator (\cite{Schm,PoPt,Argy0,HeRu,Shak}) used in the computations is

\begin{defi} \label{DEF1}
\begin{equation} \label{defi1}
[ y, x ; F]^{(1)}_{i\,j}= \left(F_i (y_1, \ldots, y_{j-1}, y_j,x_{j+1}, \ldots, x_m) - F_i (y_1, \ldots,y_{j-1}, x_j, x_{j+1}, \ldots, x_m)\right) / (y_j - x_j), \; \; \:
{\textstyle  1 \le i,j \le m .}
\end{equation}
\end{defi}

This operator defined is {\it a divided difference of $F$ on the points $x$ and $y$} because it is a bounded linear operator which satisfies the condition \cite{Schr,PoPt2}
\begin{equation} \label{cDD}
 [y,x;F] (y-x) = F(y) - F (x).
\end{equation}

\vspace{2mm}
When we evaluate $F$ in any iterative function  we calculate $m$ component functions and if we compute a divided difference then we evaluate $m(m-1)$ scalar functions, where $F(x)$ and $F(y)$ are computed separately. We must add $m^2$ quotients from any divided difference. In order to compute an inverse linear operator we have $m(m-1)(2m-1)/6$ products and $m(m-1)/2$ quotients in the decomposition $LU$ and $m(m-1)$ products and $m$ quotients in the resolution of two triangular linear systems.
Taking the previous considerations, for the first iterative function $\Phi_0$ we have $a(m) = m^2 + 2m$, $p(m)= m(2m^2+3m-5)/6+\ell\,m(3m+1)/2$ we obtain

\vspace{-2mm}
\begin{equation}\label{C0}
{\cal{C}}_0 = m(m + 2) \mu + \frac{m}{6}(2m^2+3m-5)+\frac{\ell\,m}{2}(3m+1) \quad {\rm and} \quad CEI_0 = 2^{\,1 /{\cal{C}}_0} .
\end{equation}
In an analogous way for $\Phi_1$  we get

\vspace{-2mm}
\begin{equation}\label{C1}
{\cal{C}}_1 =  2m(m + 1) \mu + \frac{m}{3}(2m^2+3m-5)+\ell\,m(3m+1) \quad {\rm and} \quad CEI_1 = 4^{\,1 /{\cal{C}}_1}.
\end{equation}
Finally, for the third iterative method we take  $a(m) = m(2m + 3)$,  $p(m)= m(2m^2+6m-8)/3+\ell\,m(3m+2)$ and the cost and the computational efficiency index are

\vspace{-2mm}
\begin{equation}\label{C2}
    {\cal{C}}_2 =  m(2m + 3) \mu + \frac{m}{3}(2m^2+6m-8)+\ell\,m(3m+2)  \quad {\rm and} \quad CEI_2 = 6^{\,1 /{\cal{C}}_2} .
\end{equation}


\vspace{2mm}

In order to compare the efficiency of the iterative methods $\Phi_i, \; 0 \le i \le 2$, we define the ratio
\begin{equation*}
R_{i,j} = \frac{\log CEI_i(\mu, m, \ell)}{\log CEI_j(\mu, m, \ell)}
  = \frac{\log(\rho_i)\: \mathcal{C}_j(\mu, m, \ell)}{\log(\rho_j)\: \mathcal{C}_i(\mu, m, \ell)}.
\end{equation*}

For $R_{i,j}>1$ the iterative method $\Phi_i$ is more efficient than $\Phi_j$. Notice that when if $R_{i,j}=1$, then we have the border between the two computational efficiencies. We study this equation for each pair of methods to obtain the following theorem:

\vspace{-1mm}
\begin{thm} \label{th1}
For all $\ell\geq 1$ and $m\ge 2$ we have:
$$
 CEI_2 > CEI_1 > CEI_0
$$
\end{thm}

\vspace{-3mm}
{\em Proof}.

{\bf 1}. In the comparison $(i,j)=(2,1)$ we have
$$
R_{\,2\,1} = \, \frac{\log \rho_2}{\log \rho_1} \: \frac{2m^2+6m\mu+3m+9\ell m+6\mu+3\ell-5}{2m^2+6m\mu+6m+9\ell m+9\mu+6\ell-8}. \\
$$
The equation $R_{\,2\,1}=1$ has one vertical asymptote for $m=0.7095$ in the $(m,\mu)$ variables and for all $\ell\ge 1$. Furthermore, $\mu$ is always negative for $m\ge2$ and for all $\ell\ge 1$. We conclude that $CEI_2 > CEI_1$ since $R_{\,2\,1}>1$ for $m\ge2$, $\mu>0$ and $\ell\ge 1$.

{\bf 2}. In the comparison $(i,j)=(1,0)$ we have
$$
R_{\,1\,0} = \frac{2m^2+6m\mu+3m+9\ell m+12\mu+3\ell-5}{2m^2+6m\mu+3m+9\ell m+6\mu+3\ell-5}. \\
$$
Equation $R_{\,1\,0}=1$ is equivalent to $6\mu=0$ for all $\ell\ge 1$ and $m\ge2$. Finally, we conclude that $CEI_1 > CEI_0$ since $R_{\,1\,0}>1$ for $m\ge2$, $\mu>0$ and $\ell\ge 1$ .
The proof is completed.
\hspace{3mm} $\Box$
%
%
%
%

\subsection{Numerical analysis}

The numerical computations were performed on MPFR library of C++ multi-precision arithmetics \cite{htt} with $4096$ digits of mantissa. All programs was compiled by {\tt gcc(4.3.3)} for i486-linux-gnu
with {\tt libgmp (v.4.2.4)} and {\tt libmpfr (v.2.4.0)} libraries on an Intel\circledR Xeon E5420, 2.5GHz and 6MB cache. For this hardware and software the computational cost of the quotient respect to the product is $\, \ell = 2.5$ (see Table \ref{TCPU_L}). Within each example the starting point is the same for all methods tested. The classical stopping criteria
$||e_I||=||x_I-\alpha|| > 0.5\cdot 10^{-\varepsilon}$ and $||e_{I+1}|| \le 0.5\cdot 10^{-\varepsilon}$,  with $\varepsilon=4096$, is replaced by
$$
E_I=\frac{||\hat{e}_I||}{||\hat{e}_{I-1}||}> 0.5\cdot 10^{-\eta} \quad {\rm and} \quad   E_{I+1} \le 0.5\cdot 10^{-\eta},
$$

where $\hat{e}_I=x_I-x_{I-1}$ and $\eta=\frac{\rho -1}{\rho^2}\: \varepsilon$. Notice that this criterium is independent of the knowledge of the root (see \cite{gng}). Furthermore, in all computations we have substituted the computational order of convergence (COC) \cite{wefe} by an approximation (ACOC) denoted by $\,\hat{\rho}$ \cite{gng} and defined as follows
$$
  \hat{\rho}=\frac{\ln E_I}{\ln E_{I-1}}\,.
$$

According to the definition of the computational cost (\ref{CEIC}), an estimation of the factor $\,\mu$ is claimed.
To do this, we express the cost of the evaluation of the elementary functions in terms of products, which depends on the machine, the software and the used arithmetics \cite{htt,fou}.
In Table \ref{TCPU_L} an estimation of the cost of the elementary functions in product units is shown, where running
time of one product is measured in milliseconds.

\vspace{3mm}
\begin{table}[htb]
\vspace{-5mm}
\caption{\small Computational cost of elementary functions computed with a program write in C++, compiled by {\tt gcc(4.3.3)} for i486-linux-gnu using {\tt libgmp (v.4.2.4)} and {\tt libmpfr (v.2.4.0)} libraries on an Intel\circledR Xeon E5420, 2.5GHz, 6MB cache, where $x = \sqrt{3}-1$ and $y = \sqrt{5}$.}
\label{TCPU_L}
\vspace{-2mm}
\begin{center}
{\footnotesize
\begin{tabular}{lcccccccc}
\hline
Arithmetics  & $x*y$ &${x}/y$ & $\sqrt{x}$& $\exp(x)$& $\ln(x)$& $\sin(x)$& $\cos(x)$& $\arctan(x)$ \\[0.1em]
\hline
C++ MPFR $4096$ digits & $0.1039$ {\small \it ms} & $2.5$    & $1.7$       & $87.8$      &  $66$    &     $116$  &    $113$   &   $228$  \\
\hline
\end{tabular}
}
\end{center}
\end{table}

\vspace{-1mm}
We analyze three examples. The first one does not present any problem when we apply the definition \ref{DEF1} given in (\ref{defi1}). The second and third examples, presented in the next section, are in fact counterexamples. That is, they need the new definition \ref{DEF2} given in (\ref{defi2}) to attain the local convergence order that theoretical results state.

\vspace{2mm}
Tables \ref{case1}--\ref{case3} show the results obtained for the iterative methods $\Phi_k$, $0\leq k\leq 2$. In each table we show the necessary iteration number $I$, the computational cost ${{\cal{C}}_k}$ in terms of products, the computational efficiency index $CEI_k$ (\ref{C0})--(\ref{C2}), the measure time factor $TF$ defined as $\: TF = {1}/{\log_{10} CEI_k}$, the necessary time $T$ in milliseconds to reach the $I$th-iteration, and the computed value of ACOC with its higher bound $\Delta \widehat{\rho_k} $, where the local order is computed by $\rho_k=\widehat{\rho_k}\pm \Delta \widehat{\rho_k}$ \cite{gng}. Furthermore, the last column shows the number of correct decimals in $x_{\scriptscriptstyle I}$, say $q$.

\vspace{2mm}
We finalize this section with a first numerical comparative study in the following system of nonlinear equations
we consider the system $F(x)=0$ with five exponential equations defined by

\vspace{-2mm}
\begin{equation}\label{exemple1}
 F(x_1,x_2,x_3,x_4,x_5)=\left\{\begin{array}{c} x_2 + x_3 + x_4 + x_5 - e^{-x_1} \,, \\
                                        x_1 + x_3 + x_4 + x_5 - e^{-x_2} \,, \\
                                        x_1 + x_2 + x_4 + x_5 - e^{-x_3} \,, \\
                                        x_1 + x_2 + x_3 + x_5 - e^{-x_4} \,, \\
                                        x_1 + x_2 + x_3 + x_4 - e^{-x_5} \,,
                        \end{array}\right.
\end{equation}
where $(m, \mu)= (5,\,87.8)$. The results shown in table \ref{case1}, where $x_{0}=(-2.1,-2.1,6.4,6.4,-2.1)^t$
 tends towards the root $\alpha\approx (-2.153967996\dots, -2.153967996\dots, 6.463463374\dots, 6.463463374\dots, -2.153967996\dots)^t$. We verify the assertion of the theorem \ref{th1}, that is, the efficiency of iterative method $\Phi_2$ is better than the efficiency of $\Phi_1$ and $\Phi_1$ one is better than $\Phi_0$ one. Moreover, we obtain the local order of convergence given in theorems \ref{th21} and \ref{th22}.

\vspace{3mm}
{\footnotesize
{\begin{table}[hh!!]
\vspace{-5mm}
 \caption{\small Numerical results for the non linear system (\ref{exemple1}).} \label{case1}
\vspace{-3mm}
 \begin{center}
 \begin{tabular}{cccccccc}
\hline \hline  \\[-0.6em]
    & $I$ & ${\cal{C}}$ & $CEI$ & $TF$ & $T$  & $ \hat{\rho}\pm\Delta\,\hat{\rho}$ & $q$ \\[0.1em]
\hline
 & \\ [-0.8em]
 $\Phi_0$ & $11$& $3223.0$ & $1.000215086$ & $10706.57$& $3730$ & $2\pm 3.04\cdot 10^{-24}$ & $3493$ \\
 $\Phi_1$ & $5$ & $5568.0$ & $1.000249006$ & $9248.26$ & $2940$ & $4\pm 7.89\cdot 10^{-11}$ & $1112$ \\
 $\Phi_2$ & $4$ & $6039.5$ & $1.000296717$ & $7761.36$ & $2540$ & $6\pm 2.00\cdot 10^{-7}$ & $1191$ \\
\hline \hline
\end{tabular}
\end{center}
\end{table}}}


\section{First-order divided difference operators preserving the local order}

This section is devoted to the preservation of the local order of convergence, of a given method,
when we replace the derivatives by divided differences. The most usual definition in practice of divided
difference is given in (\ref{defi1}).
However, we will need a new definition of divided differences instead of definition \ref{DEF1} given in  (\ref{defi1}) to get the local order wanted when we apply algorithms $\Phi_k, \; \: k=0,1,2$. We present two counterexamples
of this fact.

{\bf Case 1}. The first numerical counterexample presented consists in the following $2\times 2$ system of nonlinear equations
\vspace{-1mm}
\begin{equation} \label{exemple2}
  F (x_1, x_2) \,=\, \left\{\begin{array}{l} x^2_1 + x^2_2 -9 \,=\, 0,\\ [0.5em]  x_1\,x_2 - 1 \,=\, 0,\end{array}\right.
\end{equation}

where we have $\mu = 1.5$. When we test the behavior of the family of iterative methods $\Phi_k$, $ 0 \le k \le 2$, with initial value ${x}_{0}=(3.0,0.4)^t$ towards the root $\alpha \approx \left( 2.98118805,0.335436739 \right)^t$, the computation of ACOC give $2$, $3$ and $4$ respectively, against the result wanted. Furthermore, a comparison of the expression (\ref{equivDD}), taking into account the definition \ref{DEF1}, gives the following result
$$
\int_0^1 \: F ' ({x} + t {h}) \, dt  =  \left( \begin{array}{cc} 2 x_1 + h_1 & 2 x_2 + h_2 \\ [0.5em]
 x_2+ h_2 / 2 & x_1 + h_1 / 2 \end{array} \right) \: \mathbf{\neq} \: \left[ \, {x+h}, {x} ; F \,\right]^{(1)} = \left( \begin{array}{cc} 2 x_1 + h_1 & 2 x_2 + h_2 \\ [0.5em]
x_2  & x_1 + h_1 \end{array} \right),
$$
where $\,{x}=(x_1, x_2)^T$,  $\, {h}=(h_1, h_2)^T$ and $t \in \mathbf{R}$. Due to Potra \cite{PoPt,Pot1} we have the following necessary and sufficient condition to characterize the divided difference operator by means of a Riemann integral:
\begin{thm} If $F$ satisfies the following Lipschitz condition  $\|[x,y;F] - [u,v;F]\| \le H \left( \|x-u\| + \|y-v\|\right)$, then equality (\ref{equivDD}) holds for every pair of distinct points $(x+h,x)\in D\times D$ if and only if for all $(u,v)\in D\times D$ with $u\neq v$ and $2 v - u \in D$ it is satisfied the relation
\begin{equation} \label{cpot}
 [u,v;F] \,= \: 2\,  [u,2v-u;F] -  [v, 2v-u;F] .
 \end{equation}
\end{thm}
It can be checked that the function considered in (\ref{exemple2}) does not hold (\ref{cpot}). We need a new definition of divided differences instead of definition \ref{DEF1} given in  (\ref{defi1}) to get the local order wanted when we apply algorithms $\Phi_k, \; \: k=0,1,2$, in this case. We take the following way to compute the divided difference operator:

\vspace{2mm}
\begin{defi} \label{DEF2}
\begin{eqnarray} \nonumber
\hspace{-8mm} [ y, x ; F]^{(2)}_{i\,j} &\hspace{-4mm} =& \hspace{-4mm}\big( F_i (y_1, \ldots, y_{j-1}, y_j, x_{j+1}, \ldots, x_m)
                         - F_i (y_1, \ldots, y_{j-1}, x_j, x_{j+1}, \ldots, x_m)  \\ [0.6em] \label{defi2}
                   & & \hspace{-5mm}     + F_i (x_1, \ldots, x_{j-1}, y_j, y_{j+1}, \ldots, y_m)
                         - F_i (x_1, \ldots, x_{j-1}, x_j, y_{j+1}, \ldots, y_m)
                         \big) /\left( 2\,(y_j - x_j)\right), \; \,
{\scriptstyle  1 \le i,j \le m .}
\end{eqnarray}
\end{defi}

Notice that this operator defined in (\ref{defi2}) is not only a divided difference of $F$ on the points $x$ and $y$ because it accomplishes  the condition (\ref{cDD}) but also is a symmetric operator:  $[y,x;F] = [x,y;F]$. If we use the definition \ref{DEF2} we have to double the number of evaluations of the scalar functions in the computation of  $[y,x;F]$.

\vspace{2mm}
A development of Taylor of the expressions considered in the previous definitions gives

\vspace{-2mm}
\begin{eqnarray*}
\int_0^1 \,D_j F_i ({x} + t {h})\, dt & = & D_j F_i ({x}) + \frac{1}{2} \sum_{k=1}^m D_{k j} F_i ({x}) \, h_k + \frac{1}{6} \sum_{k,\ell=1}^m D_{k\,\ell j} F_i ({x}) \, h_k\,h_{\ell} + O({h}^3),\\ [0.2em] \label{TDDo}
 [\, {x} + {h}, {x} ; F ]^{(1)}_{i j} & = & D_j F_i ({x}) +  \sum_{k=1}^{j-1} D_{k j} F_i ({x}) \, h_k + \frac{1}{2} D_{j\,j} F_i ({x}) \, h_j + O({h}^2),\\ [0.2em]
 [\, {x} + {h}, {x} ; F ]^{(2)}_{i j} & = & D_j F_i ({x}) +  \frac{1}{2} \sum_{k=1}^m D_{k j} F_i ({x}) \, h_k + \frac{1}{4} \sum_{k,\ell=1\atop k+l<2j}^j D_{k\,\ell j} F_i ({x}) \, h_k\,h_{\ell} + \frac{1}{6} \: D_{j j j} F_i ({x}) \,h^2_j \\ [0.1em]
  & & +\: \frac{1}{4} \sum_{k,\ell=j\atop k+l>2j}^m D_{k\,\ell j} F_i ({x}) \, h_k\,h_{\ell}  + O({h}^3).
\end{eqnarray*}

From the previous expansions we obtain
$$
\int_0^1 \,D_j F_i ({x} + t {h})\, dt  - [\, {x} + {h}, {x} ; F ]^{(1)}_{i j} \,=\, O({h}),
\quad
{\rm and}
\qquad
\int_0^1 \,D_j F_i ({x} + t {h})\, dt  - [\, {x} + {h}, {x} ; F ]^{(2)}_{i j} \,=\, O({h}^2),
$$
and the definition \ref{DEF2} has higher precision. Moreover, if we use the computational definition \ref{DEF1} given in (\ref{defi1}), then the operator $\nu^{-1} = 2\,\left[ \, y, x ; F \right] -\,\left[ \, x+F(x), x-F(x) ; F \right] $, can be written as

\begin{equation*}\label{1Op}
    \nu^{-1} = \Gamma \left( 2\,I +\, 2 \widetilde{B_2} \,e + O (e^2)\right)  -\,\Gamma \Big( I  +\, 2 B_2 \,e  + O (e^2) \Big),
\end{equation*}

where $\widetilde{B_2}, B_2  \in\mathscr{L}_2\left(\mathbf{R}^m,\mathbf{R}^m\right)$, the set of bounded $2$-linear operators, but  $\widetilde{B_2}\neq A_2$ and $ \,B_2 \neq A_2$, since in (\ref{exemple2}) we have $ D_{k j} F_i (x) \neq 0$. This fact produces a loose of convergence order in many algorithms. Notice that in this case it appears the term in $e$ and the increase of local order of iterative method $\Phi_1$ (\ref{VO}) respect to $\Phi_0$ (\ref{VN}) is only one.
On the other hand, if we consider the computational definition \ref{DEF2} given in (\ref{defi2}), when we develop the operator $\nu^{-1}$,  then the difference with the theoretical result appears in the terms of second degree in $e$, while in first order are identical. Namely,

\begin{equation}\label{2Op}
    \nu^{-1} = \Gamma \left( 2\,I +\, 2 A_2 \,e + \, 2 ( \widetilde{B_3} +A^2_2 ) \,e^3 +O (e^3)\right)  -\,\Gamma \Big( I  +\, 2 A_2 \,e  + \, B_3\, ( 3 e^2 +\widetilde{e}^2 ) + O (e^3) \Big),
\end{equation}
where $\widetilde{B_3}, B_3  \in\mathscr{L}_3\left(\mathbf{R}^m,\mathbf{R}^m\right)$.

\vspace{1mm}
Using definition \ref{DEF2} for computing the divided difference operator we have to add more evaluations of functions. The new CEI (${\cal{C}}$) for $\Phi_k$ indicates by $CEI_k^{(2)}$ (${\cal{C}}_k^{(2)}$), $k=1,2$. Note that for $\Phi_0$ it is not necessary the new definition.

\vspace{-2mm}
\begin{eqnarray*}
{\cal{C}}_1^{(2)} &=& 4m^2\mu + \frac{m}{3}(2m^2+3m-5)+\ell\,m(3m+1), \quad CEI_1^{(1)} = 3^{\,1 /{\cal{C}}_1^{(1)}}, \quad CEI_1^{(2)} = 4^{\,1 /{\cal{C}}_1^{(2)}}. \\
{\cal{C}}_2^{(2)} &=& m(4m+1)\mu + \frac{m}{3}(2m^2+6m-8)+\ell\,m(3m+2), \quad CEI_2^{(1)} = 4^{\,1 /{\cal{C}}_2^{(1)}}, \quad CEI_2^{(2)} = 6^{\,1 /{\cal{C}}_2^{(2)}}.
\end{eqnarray*}

With these new computational efficiencies we obtain the following theorem:
\begin{thm} \label{th2}
For all $m\ge 2$ and $\ell\geq 1$ we have

\vspace{-3mm}
$$
 CEI_2^{(2)} > CEI_1^{(2)},
$$
but $CEI_0$ can be greater than $CEI_1^{(2)}$ and $CEI_2^{(2)}$. That is

\begin{itemize}

\vspace{-2mm}
\item For $m=2$ and for all $\ell\geq 1$ we have
$
 CEI_0 = CEI_1^{(2)} \quad \textrm{and}\quad CEI_2^{(2)} > CEI_2^{(1)}.
$

\vspace{-1mm}
\item For all $m>2$ and $\ell\geq 1$ we have
$
 CEI_0 > CEI_1^{(2)},
$

\vspace{-1mm}
\item For other comparisons it depends on the values of $\mu$, $m$ and $\ell$ (see proof).
\end{itemize}
\end{thm}

\begin{figure}[hh!!]
\begin{center}

\vspace{-7mm}
\includegraphics[width=0.60\textwidth]{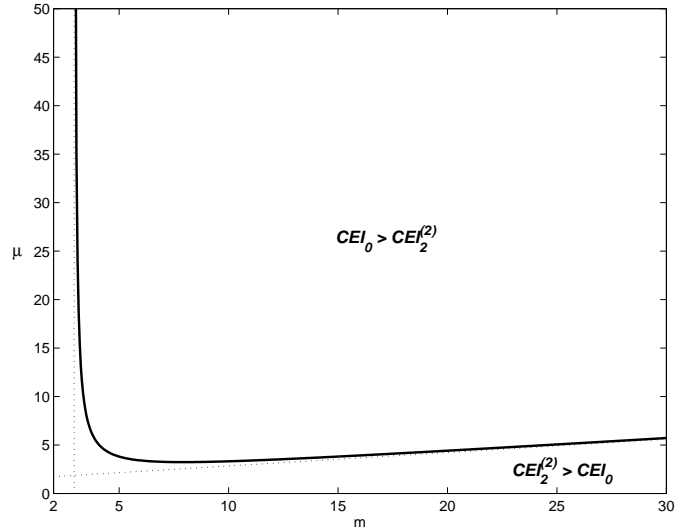}

\vspace{-5mm}
\caption{\small Boundary curve $\mu=G_{2,0}^{(2\,1)}$ for $\ell=2.5$ (\ref{cG20}).}
\label{G20}
\end{center}
\end{figure}

\vspace{-7mm}
{\em Proof}.

\vspace{1mm}
{\bf 1}. In the comparison $(i,j)=(2,1)$ we have
$$
R_{\,2\,1}^{(2\,2)} = \, \frac{\log \rho_2}{\log \rho_1} \: \frac{2m^2+12m\mu+3m+9\ell m+3\ell-5}{2m^2+12m\mu+6m+9\ell m+3\mu+6\ell-8}.
$$
The equation $R_{\,2\,1}^{(2\,2)}=1$ has one vertical asymptote for $m=0.8548$ in the $(m,\mu)$ variables and for all $\ell\ge 1$. Furthermore, $\mu$ is always negative for $m\ge2$ and for all $\ell\ge 1$. We conclude that $CEI_2^{(2)} > CEI_1^{(2)}$ since $R_{\,2\,1}>1$ for $m\ge2$, $\mu>0$ and $\ell\ge 1$.

\vspace{1mm}
{\bf 2}. In the comparison $(i,j)=(1,0)$ we have
$$
R_{\,1\,0}^{(2\,1)} = \frac{2m^2+6m\mu+3m+9\ell m+12\mu+3\ell-5}{2m^2+12m\mu+3m+9\ell m+3\ell-5}.
$$
Equation $R_{\,1\,0}^{(2\,1)}=1$ is equivalent to equation $m=2$ for all $\ell\ge 1$ and $\mu>0$. Finally, we conclude that $CEI_0 > CEI_1^{(2)}$ since $R_{\,1\,0}^{(2\,1)}<1$ for $m>2$, $\mu>0$ and $\ell\ge 1$ and $CEI_0 = CEI_1^{(2)}$ for $m=2$.

\vspace{1mm}
{\bf 3}. In the comparison $(i,j)=(2,0)$ we have
$$
R_{\,2\,0}^{(2\,1)} = \, \frac{\log \rho_2}{\log \rho_0} \: \frac{2m^2+6m\mu+3m+9\ell m+12\mu+3\ell-5}{4m^2+24m\mu+12m+18\ell m+6\mu+12\ell-16}.
$$
The particular boundary $R_{\,2\,0}^{(2\,1)}=1$ expressed by $\mu$  written as a function of $\ell$ and $m$ is
\begin{equation}\label{cG20}
G_{2,0}^{(2\,1)} = \frac{1}{3}\,\frac{2q m^2+ 3(3q\ell-r)m-3r\ell-(2q-3r)}{ 2rm-(7q+3r)},
\end{equation}
where $q=\ln(3/2)$ and $r=\ln(8/3)$. The function $G_{2,2}^{(2\,1)}$ presents a vertical asymptote for  $m = (7q+3r)/(2r) = 2.9468\ldots$ and it is always negative for $m=2$ and $\ell\ge 1$. In this case, we conclude that $CEI_2^{(2)}>CEI_0$ since $R_{\,2\,0}^{(2\,1)}>1$ for $m=2$, $\mu>0$ and $\ell\ge 1$.

Finally, for $m>2$ the boundary $\mu=G_{2,0}^{(2\,1)}$ separates two regions which extend for increasing values of $m$ and $\mu$ (see figure \ref{G20}).

\vspace{1mm}
{\bf 4}. In the comparison $(i,j)=(2,2)$ with different local order and cost we have

\vspace{-2mm}
$$
R_{\,2\,2}^{(2\,1)} = \, \frac{\log 6}{\log 4} \: \frac{2m^2+6m\mu+6m+9\ell m+9\mu+6\ell-8}{2m^2+12m\mu+6m+9\ell m+3\mu+6\ell-8}.
$$
Here the boundary $R_{\,2\,2}^{(2\,1)}=1$ expressed by $\mu$ as a function of $\ell$ and $m$ can be written as

\vspace{-2mm}
\begin{equation}\label{cG22}
G_{2,2}^{(2\,1)} = \frac{1}{3}\,\frac{2q m^2+ 3q(3\ell+2)m+6q\ell-8q}{ 2rm-(5q+2r)}.
\end{equation}
This function, $\mu=G_{2,2}^{(2\,1)}$, has a vertical asymptote for  $m = (5q+2r)/(2r) = 2.0334\ldots$ and it is always negative for $m=2$ and $\ell\ge 1$. In this case, we conclude that $CEI_2^{(2)}>CEI_2^{(1)}$ since $R_{\,2\,2}^{(2\,1)}>1$ for $m=2$, $\mu>0$ and $\ell\ge 1$.

\begin{figure}[hh!!]
\begin{center}

\vspace{-3mm}
\includegraphics[width=0.60\textwidth]{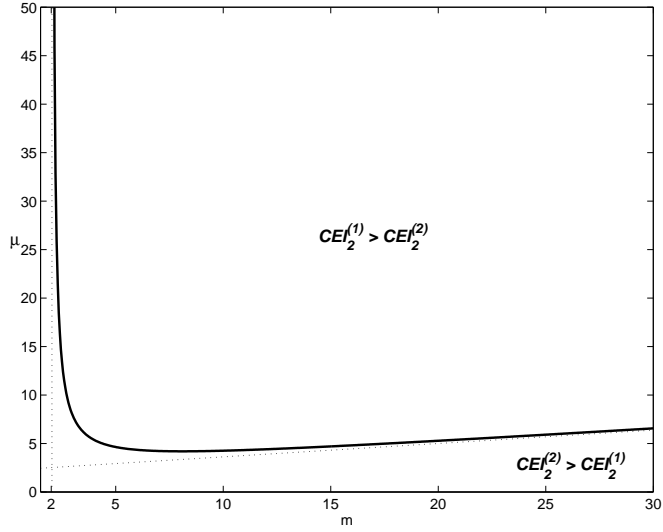}

\vspace{-5mm}
\caption{\small Boundary curve $\mu=G_{2,2}^{(2\,1)}$ for $\ell=2.5$ (\ref{cG22}).}
\label{G22}
\end{center}
\end{figure}

Finally, for $m>2$ the boundary $\mu=G_{2,2}^{(2\,1)}$ separates two regions which extend for increasing values of $m$ and $\mu$ (see figure \ref{G22}).

\vspace{1mm}
{\bf 5}. In the comparison $(i,j)=(1,1)$ with different local order and cost we have
$$
R_{\,1\,1}^{(2\,1)} = \, \frac{\log 4}{\log 3} \: \frac{2m^2+6m\mu+3m+9\ell m+6\mu+3\ell-5}{2m^2+12m\mu+3m+9\ell m+3\ell-5}.
$$
The particular boundary $R_{\,1\,1}^{(2\,1)}=1$ expressed by $\mu$  written as a function of $\ell$ and $m$ is

\begin{figure}[hh!!]
\begin{center}

\vspace{-5mm}
\includegraphics[width=0.60\textwidth]{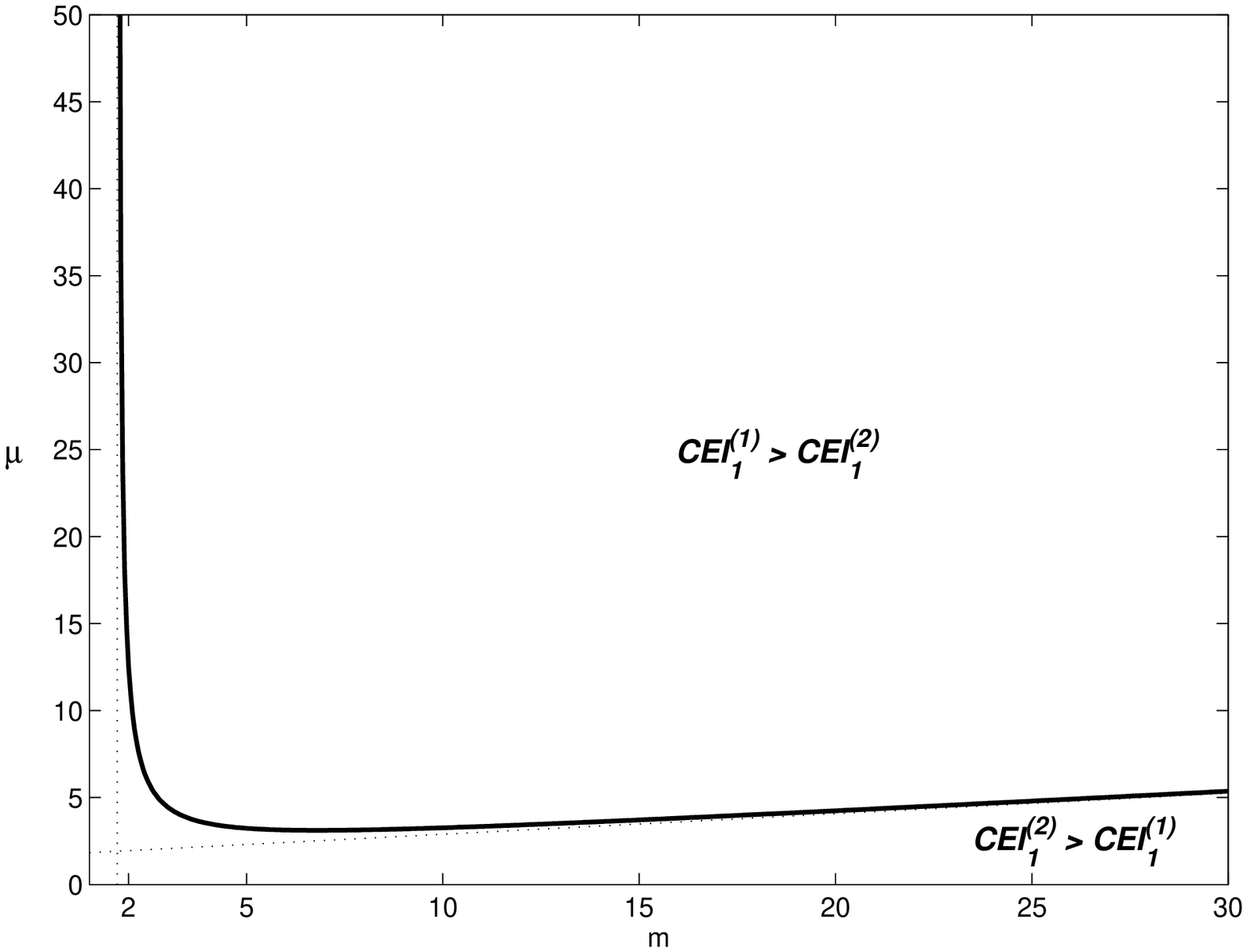}

\vspace{-5mm}
\caption{\small Boundary curve $\mu=G_{1,1}^{(2\,1)}$ for $\ell=2.5$ (\ref{cG11}).}
\label{G11}
\end{center}
\end{figure}

\vspace{-2mm}
\begin{equation}\label{cG11}
G_{1,1}^{(2\,1)} = \frac{1}{12}\,\frac{2s m^2+ 3s(3\ell+1)m+3s\ell-5s}{(t-s)m-t},
\end{equation}
where $s=\ln(4/3)$ and $t=\ln 2$.The function $G_{1,1}^{(2\,1)}$ has a vertical asymptote for  $m = t/(t-s) = 1.7095\ldots$. The boundary $\mu=G_{1,1}^{(2\,1)}$ separates two regions that are shown in figure \ref{G11}.

\vspace{1mm}
The proof is complete.  \hspace{3mm} $\Box$


{\footnotesize
{\begin{table}[hh!!]
\vspace{-7mm}
 \caption{\small Numerical results for the non linear system (\ref{exemple2}).} \label{case2}

\vspace{-2mm}
 \begin{center}
 \begin{tabular}{cccccccc}
\hline \hline  \\[-0.6em]
    & $I$ & ${\cal{C}}$ & $CEI$ & $TF$ & $T$  & $ \hat{\rho}\pm\Delta\,\hat{\rho}$ & $q$ \\[0.1em]
\hline
 & \\ [-0.8em]
 $\Phi_0$       & $11$& $32.5$ & $1.021556664$ & $107.96$ & $40$ & $2\pm 2.61\cdot 10^{-04}$ & $3334$ \\
 $\Phi_1^{(1)}$ & $7$ & $59.0$ & $1.018794991$ & $123.66$ & $40$ & $3\pm 4.07\cdot 10^{-40}$ & $2908$ \\
 $\Phi_1^{(2)}$ & $5$ & $65.0$ & $1.021556664$ & $107.96$ & $30$ & $4\pm 5.86\cdot 10^{-13}$ & $1951$ \\
 $\Phi_2^{(1)}$ & $5$ & $69.0$ & $1.020294410$ & $114.61$ & $35$ & $4\pm 1.70\cdot 10^{-23}$ & $1384$ \\
 $\Phi_2^{(2)}$ & $4$ & $75.0$ & $1.024177781$ & $ 96.38$ & $20$ & $6\pm 1.90\cdot 10^{-08}$ & $2392$ \\
\hline \hline
\end{tabular}
\end{center}
\end{table}}}


\vspace{2mm}
Table \ref{case2} shows the same columns that table \ref{case1} but taking into account the two definitions of divided difference operator. Note that for this case 2, we have $(m,\mu,\ell)=(2,1.5,2.5)$ and the theorem \ref{th2} is accomplished. That is,

\vspace{-2mm}
$$
 CEI_2^{(2)}>CEI_1^{(2)}=CEI_0>CEI_1^{(1)} \quad {\rm and} \quad CEI_2^{(2)}>CEI_2^{(1)}.
$$

Observe that the values of ACOC in table \ref{case2} confirm the theoretical local orders of convergence with a high precision.

\vspace{1mm}
{\bf Case 2}.  In this case we consider an example where the computational definition of the divided difference operator presents non zero values for $\widetilde{B_3}$ and $B_3$ (see (\ref{2Op})), opposite to the case 1.

\vspace{-1mm}
\begin{equation}\label{exemple3}
       x_i-\cos\left( 2\,x_i-\sum_{j=1}^{m} \: x_j\right) =0,  \; \quad \;  1\le i \le m,
\end{equation}

where $m= 3 $ and $\mu =  113.3 $. We test the behavior of the iterative methods $\Phi_k$, $ 0 \le k \le 2$, say (\ref{VN}), (\ref{VO}) and (\ref{VO2}),  with the initial
value  $x_{0}=(0.4,0.4,0.9)^t$ towards the root $\alpha \approx \left( 0.5438500415,0.5438500415,0.9957781534 \right)^t$, considering the two definitions of the divided difference operator presented in this paper.

{\footnotesize
{\begin{table}[hh!!]
\vspace{-2mm}
 \caption{\small Numerical results for the non linear system (\ref{exemple3}).} \label{case3}
\vspace{-3mm}
 \begin{center}
 \begin{tabular}{cccccccc}
\hline \hline  \\[-0.6em]
    & $I$ & ${\cal{C}}$ & $CEI$ & $TF$ & $T$  & $ \hat{\rho}\pm\Delta\,\hat{\rho}$ & $q$ \\[0.1em]
\hline
 & \\ [-0.8em]
 $\Phi_0$       & $13$& $1748.0$ & $1.000396616$ & $5806.73$ & $2330$ & $2\pm 1.58\cdot 10^{-15}$ & $2575$ \\
 $\Phi_1^{(1)}$ & $8$ & $2816.2$ & $1.000390181$ & $5902.48$ & $2310$ & $3\pm 9.56\cdot 10^{-35}$ & $2549$ \\
 $\Phi_1^{(2)}$ & $6$ & $4175.8$ & $1.000332038$ & $6935.85$ & $2580$ & $4\pm 5.88\cdot 10^{-07}$ & $2517$ \\
 $\Phi_2^{(1)}$ & $6$ & $3169.6$ & $1.000437468$ & $5264.59$ & $2000$ & $4\pm 5.35\cdot 10^{-08}$ & $1514$ \\
 $\Phi_2^{(2)}$ & $4$ & $4529.2$ & $1.000395680$ & $5820.46$ & $1870$ & $6\pm 7.52\cdot 10^{-06}$ & $725$ \\
\hline \hline
\end{tabular}
\end{center}
\end{table}}}


Table \ref{case3} shows the same columns that table \ref{case2}. Notice that the values of ACOC confirm the theoretical local orders of convergence of the iterative methods presented with a high precision.  In this case 3, we have $(m,\mu,\ell)=(3,113.3,2.5)$ and the theorem \ref{th2} is accomplished. That is,

\vspace{-2mm}
$$
 CEI_2^{(1)}>CEI_0>CEI_2^{(2)}>CEI_1^{(2)} \quad {\rm and} \quad CEI_1^{(1)}>CEI_1^{(2)}.
$$

\section{Concluding remarks}

From Ostrowski's method witch is a root finder for nonlinear equations we present two modified algorithms derivative-free for solving systems of nonlinear equations.
A development of an inverse first-order divided difference operator for functions of several variables is analyzed.
Moreover, we study the computational efficiency index and a comparison for the three iterative methods is given.
A similar analysis should be performed for other Newton's type methods obtaining new derivative-free algorithms.

\vspace{2mm}
Numerical examples that illustrate the theoretical results presented in this paper are also given.
We emphasize the necessity sometimes of a new way to compute the matrix of divided difference operator since in many examples the theoretical one and the classical way used until now to calculate it do not coincide. With the new definition the estimate of the computational efficiency index is revisited. Of course, other possible examples
of divided differences should be considered.

\vspace{2mm}
An approximation of the root with high precision is obtained if we use a multi-precision arithmetics.
To illustrate the technique presented three examples are worked an completely solved.
An approximation of the computational order of convergence is computed independently of the knowledge of the root and it confirms in a numerical way that the convergence order of the methods is well deduced.

\section*{Acknowledgements}

Research supported by MTM2011-28636-C02-01 project of the Spanish Ministry of
Science and Innovation for the first two authors. The third author's research has supported by
Spanish grants MTM2007-62945 and 08662/PI/08.


\end{document}